\documentclass[12pt,a4paper]{amsart}
\usepackage{t1enc}
\usepackage[latin1]{inputenc}
\usepackage[english]{babel}
\usepackage{hyperref}
\usepackage{graphicx}
\usepackage{latexsym}
\usepackage{amsmath,amsfonts,amssymb,amsthm}
\usepackage{mathptmx} 
\begin{document}
\title[Maximal Plurisubharmonic Models]{Maximal Plurisubharmonic Models}
\author{Giuseppe Tomassini}
\address{G. Tomassini: Scuola Normale Superiore, Piazza dei Cavalieri, 7 --- I-56126 Pisa, Italy}
\email{g.tomassini@sns.it}
\author{Sergio Venturini}
\address{S. Venturini: Dipartimento Di Matematica, Universit\`{a} di Bologna, \,\,Piazza di Porta S. Donato 5 ---
I-40127 Bologna, Italy}
\email{venturin@dm.unibo.it}
\keywords{Complex manifolds $\cdot$ plurisubharmonic functions $\cdot$ complex Monge-Amp\`ere equation $\cdot$ Finsler geometry}
\subjclass[2000]{Primary 32F45, 32C09 Secondary 32Q45, 32T35}
{\large 


%
%
\def\R{{\rm I\kern-.185em R}}
\def\C{{\rm\kern.37em\vrule height1.4ex width.05em depth-.011em\kern-.37em C}}
\def\N{{\rm I\kern-.185em N}}
\def\Z{{\bf Z}}
\def\Q{{\mathchoice{\hbox{\rm\kern.37em\vrule height1.4ex width.05em 
depth-.011em\kern-.37em Q}}{\hbox{\rm\kern.37em\vrule height1.4ex width.05em 
depth-.011em\kern-.37em Q}}{\hbox{\sevenrm\kern.37em\vrule height1.3ex 
width.05em depth-.02em\kern-.3em Q}}{\hbox{\sevenrm\kern.37em\vrule height1.3ex
width.05em depth-.02em\kern-.3em Q}}}}
\def\P{{\rm I\kern-.185em P}}
\def\H{{\rm I\kern-.185em H}}
%
\def\Aleph{\aleph_0}
\def\ALEPH#1{\aleph_{#1}}
\def\sset{\subset}\def\ssset{\sset\sset}
%
\def\bar#1{\overline{#1}}
\def\dim{\mathop{\rm dim}\nolimits}
\def\half{\textstyle{1\over2}}
\def\Half{\displaystyle{1\over2}}
\def\mlog{\mathop{\half\log}\nolimits}
\def\Mlog{\mathop{\Half\log}\nolimits}
\def\Det{\mathop{\rm Det}\nolimits}
\def\Hol{\mathop{\rm Hol}\nolimits}
\def\Aut{\mathop{\rm Aut}\nolimits}
\def\Re{\mathop{\rm Re}\nolimits}
\def\Im{\mathop{\rm Im}\nolimits}
\def\Ker{\mathop{\rm Ker}\nolimits}
\def\Fix{\mathop{\rm Fix}\nolimits}
\def\Res{\mathop{\rm Res}\nolimits}
\def\sp{\mathop{\rm sp}\nolimits}
\def\id{\mathop{\rm id}\nolimits}
\def\Trace{\mathop{\rm Tr}\nolimits}
\def\cancel#1#2{\ooalign{$\hfil#1/\hfil$\crcr$#1#2$}}
\def\prevoid{\mathrel{\scriptstyle\bigcirc}}
\def\void{\mathord{\mathpalette\cancel{\mathrel{\scriptstyle\bigcirc}}}}
\def\n{{}|{}\!{}|{}\!{}|{}}
\def\abs#1{\left|#1\right|}
\def\norm#1{\left|\!\left|#1\right|\!\right|}
\def\nnorm#1{\left|\!\left|\!\left|#1\right|\!\right|\!\right|}
%
\def\upperint{\int^{{\displaystyle{}^*}}}
\def\lowerint{\int_{{\displaystyle{}_*}}}
\def\Upperint#1#2{\int_{#1}^{{\displaystyle{}^*}#2}}
\def\Lowerint#1#2{\int_{{\displaystyle{}_*}#1}^{#2}}
%
\def\rem #1::#2\par{\medbreak\noindent{\bf #1}\ #2\medbreak}
\def\proclaim #1::#2\par{\removelastskip\medskip\goodbreak{\bf#1:}
\ {\sl#2}\medskip\goodbreak}
\def\ass#1{{\rm(\rmnum#1)}}
\def\assertion #1:{\Acapo\llap{$(\rmnum#1)$}$\,$}
\def\Assertion #1:{\Acapo\llap{(#1)$\,$}}
\def\acapo{\hfill\break\noindent}
\def\Acapo{\hfill\break\indent}
\def\proof{\removelastskip\par\medskip\goodbreak\noindent{\it Proof.\/\ }}
\def\prova{\removelastskip\par\medskip\goodbreak
\noindent{\it Dimostrazione.\/\ }}
\def\qed{{\bf //}\par\smallskip}
\def\BeginItalic#1{\removelastskip\par\medskip\goodbreak
\noindent{\it #1.\/\ }}
\def\iff{if, and only if,\ }
\def\sse{se, e solo se,\ }
\def\rmnum#1{\romannumeral#1{}}
\def\Rmnum#1{\uppercase\expandafter{\romannumeral#1}{}}
\def\smallfrac#1/#2{\leavevmode\kern.1em
\raise.5ex\hbox{\the\scriptfont0 #1}\kern-.1em
/\kern-.15em\lower.25ex\hbox{\the\scriptfont0 #2}}
%
\def\Left#1{\left#1\left.}
\def\Right#1{\right.^{\llap{\sevenrm
\phantom{*}}}_{\llap{\sevenrm\phantom{*}}}\right#1}
%
\def\a{\alpha}
\def\bg{\beta}
\def\g{\gamma}
\def\G{\Gamma}
\def\dg{\delta}
\def\D{\Delta}
\def\e{\varepsilon}
\def\eps{\epsilon}
\def\z{\zeta}
\def\th{\theta}
\def\T{\Theta}
\def\k{\kappa}
\def\lg{\lambda}
\def\Lg{\Lambda}
\def\m{\mu}
\def\n{\nu}
\def\r{\rho}
\def\s{\sigma}
\def\Sg{\Sigma}
\def\ph{\varphi}
\def\Ph{\Phi}
\def\x{\xi}
\def\om{\omega}
\def\Om{\Omega}


\newtheorem{theorem}{Theorem}[section]
\newtheorem{proposition}{Proposition}[section]
\newtheorem{lemma}{Lemma}[section]
\newtheorem{corollary}{Corollary}[section]
\newtheorem{remark}{Remark}[section]
\newtheorem{definition}{Definition}[section]

\newtheorem{teorema}{Teorema}[section]
\newtheorem{proposizione}{Proposizione}[section]
\newtheorem{corollario}{Corollario}[section]
\newtheorem{osservazione}{Osservazione}[section]
\newtheorem{definizione}{Definizione}[section]
\newtheorem{esempio}{Esempio}[section]
\newtheorem{esercizio}{Esercizio}[section]
\newtheorem{congettura}{Congettura}[section]

\bibliographystyle{abbrv}
\def\references{\bibliography{bibfile}}
\def\Levi{{\mathcal L}}
\def\Ampere{Amp\`ere}
\def\APair[#1,#2]{({#1},{#2})}
\def\PSHU[#1,#2]{{\mathcal U}({#1},{#2})}
\def\CManA{M}
\def\CManB{N}
\def\SManA{V}
\def\SManB{W}
\def\mink{\mu}
\def\Cvx{D}
\def\PointA{p}
\def\PointB{q}
\def\TanVA{\xi}
\def\TanSpace#1{T{#1}}
\def\TanSpaceP[#1,#2]{T_{#1}{#2}}
\def\arctanh{\mathop{\rm arctanh}\nolimits}

\def\statementHopf{proposition}
\def\statementSchwarz{theorem}

\nocite{article:GTV}
\nocite{article:BedfordMASurvey}
\nocite{article:Burns1982A}
\nocite{article:BurnsAndHind}
\nocite{book:KobayashiHyperbolicComplexSpace}
\nocite{article:LempertEllipticTube}
\nocite{article:LempertSzoke}
\nocite{article:AbatePatrizioConstCurv}
\nocite{article:PatrizioWong}
\nocite{article:Szoke}
\nocite{book:BerensteinGay}
 
\begin{abstract}
An \emph{analytic pair} of dimension $n$ and {\it center $V$} is a pair $(\SManA,\CManA)$ where $\CManA$ is a complex manifold of (complex) dimension $n$ and $\SManA\sset\CManA$ is a closed totally real analytic submanifold of dimension $n$. To an analytic pair $(\SManA,\CManA)$ we associate the class $\PSHU[\SManA,\CManA]$ of the functions $u:\CManA\to[0,\pi/4[$ which are plurisubharmonic in $\CManA$ and such that $u(p)=0$ for each $p\in\SManA$. If $\PSHU[\SManA,\CManA]$ admits a maximal function $u$, the triple $(\SManA,\CManA,u)$ is said to be a {\it maximal plurisubharmonic model}. After defining a pseudo-metric $E_{\SManA,\CManA}$ on the center $V$ of an analytic pair $(\SManA,\CManA)$ we prove (see Theorem \ref{thm::MAMaximal}, Theorem \ref {LI}) that maximal plurisubharmonic models provide a natural generalization of the Monge-Amp\`ere models introduced by Lempert and Sz\"oke in \cite{article:LempertSzoke}.
\end{abstract}
\maketitle

\section{\label{section:Introduzione}Introduction}
An \emph{analytic pair} of dimension $n$ is a pair $(\SManA,\CManA)$ where $\CManA$ is a complex manifold of (complex) dimension $n$ and
$\SManA\sset\CManA$ is a closed totally real analytic submanifold of dimension $n$. The submanifold $\SManA$ is said to be the \emph{center} of the analytic pair $(\CManA,\SManA)$. 
 We denote by $T\CManA$, $T\SManA\sset T\CManA$ 
the respective (real) tangent fibre bundles and $J:T\CManA\to T\CManA$ the complex structure of $\CManA$. 

To an analytic pair $(\SManA,\CManA)$ we associate the class $\PSHU[\SManA,\CManA]$ of the functions
$$
u:\CManA\to[0,\pi/4[
$$
which are plurisubharmonic in $\CManA$ and such that $u(p)=0$ for each $p\in\SManA$. The choice of the constant $\pi/4$ will be explained later.

A function $u\in\PSHU[\SManA,\CManA]$ is said to be \emph{maximal} (for the pair $(\SManA,\CManA)$) if
$$
v(p)\leq u(p)
$$ 
for every $p\in\CManA$, $v\in\PSHU[\SManA,\CManA]$,
$u$ vanishes exactly on $\SManA$, that is $u(p)>0$ for every $p\in\CManA\setminus\SManA$, and
$$
\sup_{p\in\CManA}u(p)=\frac{\pi}{4}.
$$
Clearly, a maximal element in $u\in\PSHU[\SManA,\CManA]$ is unique, provided it exists.
\noindent
We say that a triple $(\SManA,\CManA,u)$ is a \emph{maximal plurisubharmonic model} (of bounded type), for short a \emph{maximal model}, if $(\SManA,\CManA)$ is an analytic pair and $u\in\PSHU[\SManA,\CManA]$ a maximal function.

With a (little) abuse of language we say that an analytic pair $(\SManA, \CManA)$ is a \emph{bounded maximal model} provided there exists a maximal function $u\in\PSHU[\SManA,\CManA]$.

Let now $u\in\PSHU[\SManA,\CManA]$ where $(V,M)$ is an analytic pair. For $\PointA\in\SManA$ and $\TanVA\in\TanSpace\SManA\sset\TanSpace\CManA$ the formula
$$
	E_{u,\CManA}(\PointA,\TanVA)=
		{\rm ``slope\ of\ }u\ {\rm at\ }\PointA\ {\rm in\ the\ direction\ }J\TanVA\ ''
$$
defines a pseudo-metric on $\SManA$ associated to the function $u$.

Taking 
$$
\sup_{u\in\PSHU[\SManA,\CManA]}E_{u,\CManA}(\PointA,\TanVA)
$$
we define a pseudo-metric $E_{\SManA,\CManA}$ on $\SManA$
which depends only on the geometry of $\PSHU[\SManA,\CManA]$. If $(\SManA,\CManA,u)$ is a maximal model $E_{\SManA,\CManA}$ actually coincides with
$E_{u,\CManA}$. (See Section \ref{section:LaMetrica}) for the precise definitions).

We now explain the motivations of our contruction.

Following \cite{article:LempertSzoke} we recall that a
\emph{Monge-Amp\`ere model of dimension $n$} is a triple
$(\SManA,\CManA,u)$ where
\begin{itemize} 
\item[1)] $(\SManA,\CManA)$ is an analytic pair of dimension $n$;
\item[2)] $u$ is a continuous, plurisubharmonic function such that $\SManA=\{u=0\}$;
\item[3)] $u$ is a smooth solution on $\CManA\smallsetminus\SManA$ of the (complex) Monge-Amp\`ere equation
$$
(dd^cu)^n=0;
$$
\item[4)] $u^2$ is smooth and strictly plurisubharmonic exhaustion function on $\CManA$.
\end{itemize}
In such conditions $V$ is called the {\it center} of the Monge-Amp\`ere model $(\SManA,\CManA,u)$. Moreover, if the function $u$ is bounded then $(\SManA,\CManA,u)$ is said to be of {\it bounded type}.

A {\it holomorphic map} $F:(\SManA_1,\CManA_1,u_1)\to(\SManA_2,\CManA_2,u_2)$ between Monge-Amp\`ere models is a holomorphic map $F:\CManA_1\to\CManA_2$ such that $F(V_1)\subset V_2$ and $u_1=F\circ u_2$.

Two Monge-Amp\`ere models
$(\SManA_1,\CManA_1,u_1)$ and $(\SManA_2,\CManA_2,u_2)$
are said to be {\it isomorphic} if there exists a biholomorphic map
$F:\CManA_1\to\CManA_2$ such that $F(\SManA_1)=\SManA_2$ and $u_1=F\circ u_2$. 

The center $V$ of a Monge-Amp\`ere model $(\SManA,\CManA,u)$ is a Riemannian manifold with metric $g$ given by the restriction to the tangent bundle $T\SManA$ of the Levi form $\mathcal L(u^2)$ of $u^2$.

Lempert and Sz\"oke proved in \cite{article:LempertSzoke} that every compact Riemannian manifold $(\SManA,g)$ is, canonically, the center of a Monge-Amp\`ere model of bounded type $(\SManA,\CManA,u)$. Moreover, $(\SManA,\CManA,u)$ is completely determined (up to isomorphisms) by the Riemannian manifold $(\SManA,g)$, i.e. two bounded Monge-Amp\`ere models $(\SManA_1,\CManA_1,u_1)$ and $(\SManA_2,\CManA_2,u_2)$ are isomorphic if and only if their respective centers $(\SManA_1,g_1)$ and $(\SManA_2,g_2)$ are isometric Riemannian manifolds . 

The canonical model is constructed as follows. Let  $u$ be the length function $\abs{\phantom{x}}:T\SManA\to[0,+\infty[$, associated to $g$.
Identify $\SManA$ with the zero section of $T\SManA$, and consider, for $0<r\leq+\infty$, the $r$-{\it tube} 
$$
\mathcal T_r\SManA=\bigl\{\TanVA\in T\SManA\mid u(\TanVA)<r\bigr\}
$$
with center $\SManA$. Then, for $r>0$ small enough, $\mathcal T_r\SManA$ carries an unique complex structure such that the triple
$(\SManA,\mathcal T_r\SManA,u)$ is a Monge-Amp\`ere model and
the restriction to the tangent bundle $T\SManA$
of the form $2\mathcal L(u^2)$ is exactly the Riemannian metric $g$ (see \cite{article:LempertSzoke} and \cite{article:Szoke},
or \cite{article:GuilleminStenzelI}).

The manifold $\mathcal T_r\SManA$ is called a \emph{Grauert tube of radius $r$} over the Riemannian manifold $\SManA$. The name ``Grauert tube'' is due to the following theorem proved by Grauert in \cite{article:GrauertTube}: every real analytic manifold $V$ of dimension $n$ 
embeds as a maximal totally real submanifold of an $n$-dimensional complex manifold $M$ in such a way to have a basis of Stein neighbourhoods.

A Grauert tube $\mathcal T_r\SManA$ is said to be \emph{rigid} if each biholomorphic automorphism
$f:\mathcal T_r\SManA\to\mathcal T_r\SManA$ preserves the center $\SManA$.

Grauert tubes, and their extension to non compact centers, are widely studied complex manifolds,
expecially in connection with curvature problems (\cite{article:LempertSzoke}) and rigidity problems
(see e. g. \cite{article:BurnsAndHind},  \cite{article:BurnsAndHalvAndHind},
\cite{article:KanSJ_A}, \cite{article:KanSJAndMa} and \cite{article:KanSJRigidA}).

By the way it would be interesting to have an analogous of the canonical model starting from a center equipped with a Finsler metric.  

The goal of this paper is to show that maximal models of bounded type provide a natural generalization of the (bounded) Monge-Amp\`ere ones. The results obtained here must be considered as a preliminary exploration of the geometry of such models. Clearly, assuming no kind of regularity of $u^2$ the Riemannian geometry (of the center) should be replaced by a ''pseudo-metric geometry''. We claim that the pseudo-metric $E_{\SManA,\CManA}$ defined in this paper is the right object for our scope.

The paper is organized as follows.

In Section \ref{section:Subarmoniche}, for the sake of completeness, we prove some
simple variations of Hopf lemma and Phragmen-Lindel\"of principle for subharmonic function
of one complex variable, in a form that we need in the sequel.

In Section \ref{section:LaMetrica} we introduce the pseudo-metrics $E_{u,\CManA}$ and
$E_{\SManA,\CManA}$ and describe their basic properties.
It turns out that if $\CManA$ is the unit disc $\Delta=\{z\in\C\mid\abs{z}<1\}$
and $\SManA=]-1,1[$, then the associated metric on the center $]-1,1[$
is the restriction of the Poincar\'e metric on $\Delta$
(this is the reason for the constant $\pi/4$ above).
Moreover, if $(\SManA_1,\CManA_1)$, $(\SManA_2,\CManA_2)$ are analytic pairs and  $F:\CManA_1\to\CManA_2$ is a holomorphic map which such that $F(V_1)\subset V_2$ then $F$ is a contraction for the corresponding pseudo-metrics on the centers.
Thus, our theory is a Kobayashi-like pseudo-metric theory.
In \cite{article:GTV} it was proved that
the class of all Finsler pseudo-metrics 
on the center of an analytic pair $\APair[\SManA,\CManA]$
having this contraction property admits a largest element $F_{\SManA,\CManA}$,
so that $E_{\SManA,\CManA}\leq F_{\SManA,\CManA}$.
For the definition and the main properties of the metric $F_{\SManA,\CManA}$ we refer to \cite{article:GTV}.
It turns out that the equality $E_{\SManA,\CManA}=F_{\SManA,\CManA}$
is related to the existence of ``complex geodesic'' for such pseudo-metrics (see Theorem \ref{thm::EGeodesics}).

It should be observed that the pseudo-metric $E_{\SManA,\CManA}$ is positively homogeneous but in general it is not symmetric, that is, for $p\in\SManA$ and $\xi\in T_p\SManA$ it may happen that $E_{\SManA,\CManA}(p,-\xi)\neq E_{\SManA,\CManA}(p,\xi)$.

Sections \ref{section:MongeAmpere}, Section \ref{section:CasoLiscio} are devoted to the interplay between maximal functions for analytic pairs $(V,M)$ and solutions of the complex Monge-\Ampere\ equation on $M\setminus V$.
After proving that if $u\in\PSHU[\SManA,\CManA]$ is a continuous exhaustion function on $\CManA$ then $(\SManA,\CManA,u)$ is a maximal model if and only if $(dd^cu)^n=0$ on $M\setminus V$ (see Theorem \ref{thm::MAMaximal}), in Section \ref{section:CasoLiscio} we show that for a Monge-Amp\`ere model $(\SManA,\CManA,u)$ of bounded type one has 
$$
E_{\SManA,\CManA}(p,\TanVA)=F_{\SManA,\CManA}(p,\TanVA)=\sqrt{2\Levi_{u^2}(p,\TanVA)}.
$$
i.e. the pseudo-metric $E_{\SManA,\CManA}(p,\TanVA)$ coincides with the Riemannian metric on $\SManA$ (see Theorem \ref {LI}).

Finally, in the last two sections we give two significnt examples of generalized Monge-Amp\`ere models maximal model which are not (unless exceptional cases) Monge-Amp\`ere models.

In {Section\llap{ }} \ref{section:esa}
we prove that if $\mu:\R^n\to[0,+\infty[$
is the Minkowsky funcional associated to a bounded open
convex subset of $\R^n$ containing the origin
(not necessarily symmetric with respect to the origin)
then $(\R^n,X_\mink,u)$, where
$$
X_\mink=\{z=x+iy\in\C^n\mid\mink(y)<\pi/2\} 
$$
and $u(z)=u(x+iy)=\mink(y)$, is a maximal model.
It is worthy of observing that  this example easily generalizes  
if $\R^n$ is replaced by an arbitrary real Banach space where, in general,  we have no analogous of the Monge-Amp\`ere operator, while the definition of maximal plurisubharmonic bounded function is exactly the same.

In Section \ref{section:TuboEllittico} we shall prove that
if $\Cvx\sset\R^n$ is a bounded open convex set and $\Cvx^{\rm ell}$  the elliptic tube over the convex set $\Cvx$ described by Lempert in \cite{article:LempertEllipticTube}, then $(\Cvx,\Cvx^{\rm ell})$ is a bounded maximal model for which an explicit description of the corresponding maximal function $u$ is provided.
\section{\label{section:Subarmoniche}Subharmonic functions}
\def\MyStripAB{S}
\def\MyDomain{D}
\def\UpperStrip{r}
\def\UpperFunc{a}

The pourpose of this section is to prove following Theorem \ref{thm::PSHSchwarz}, which is the
``Schwarz lemma'' in our context.

Let us begin recalling the classical Hopf lemma
in the form that we need in the sequel.
\begin{\statementHopf}\label{thm::SHHopfLemmaA}
Let  $\MyDomain\sset\C$ be open, $\MyDomain\neq\C$, and $u:\MyDomain\to[-\infty,0[$
a negative subharmonic function.
For $z\in\MyDomain$ denote by $\delta(z)$ the distance from z to $\partial\MyDomain$.
Let $x\in\partial\MyDomain$ and assume that
$\partial\MyDomain$ is of class $C^2$ in a neighbourhood of $x$. Then
$$
	\limsup_{\MyDomain\ni z\to x}\frac{u(z)}{\delta(z)}<0.
$$
\end{\statementHopf}

For a proof see e.g. Proposition 12.2 of \cite{article:FornaessStensones}.

\begin{\statementSchwarz}\label{thm::PSHSchwarz}
Given $r,a>0$,
let
$$
\MyDomain=\bigl\{z\in\C\mid0<\Im z<\UpperStrip\bigr\}
$$
and $u:\MyDomain\to[0,a[$ a bounded subharmonic function such that
$$
	\lim_{\MyDomain\ni z\to x}u(z)=0
$$
for each $x\in\R$.
Then, for every $z\in\MyDomain$, $x\in\R$
$$
u(z)\leq\frac{\UpperFunc}{\UpperStrip}\Im z, \>\> \limsup_{y\to0^+}\frac{u(x+iy)}{y}\leq\frac{\UpperFunc}{\UpperStrip}.
$$
If there exist either 
$z_0\in\MyDomain$ such that
$$
u(z_0)=\frac{\UpperFunc}{\UpperStrip}\Im z_0
$$
or $x_0\in\R$ such that
$$
\limsup_{y\to0^+}\frac{u(x_0+iy)}{y}=\frac{\UpperFunc}{\UpperStrip},
$$
then
$$
u(z)=\frac{\UpperFunc}{\UpperStrip}\Im z
$$
for every $z\in\MyDomain$.
\end{\statementSchwarz}

\proof
The function $v:\MyDomain\to\R$ defined by
$$
v(z)=u(z)-\frac{\UpperFunc}{\UpperStrip}\Im z
$$
is bounded, subharmonic on $\MyDomain$ and satisfy
$$
\limsup_{\MyDomain\ni z\to\z}v(z)\leq0
$$
for each $\z\in\partial\MyDomain$.

By the Phragmen-Lindel\"of principle
(see e.g. Proposition 4.9.45, pag. 463 of \cite{book:BerensteinGay})
for each $z\in\MyDomain$
$$
v(z)\leq0,
$$
that is
$$
	u(z)\leq\frac{\UpperFunc}{\UpperStrip}\Im z,
$$
and hence, for each $x\in\R$,
$$
	\limsup_{y\to0^+}\frac{u(x+iy)}{y}=\frac{\UpperFunc}{\UpperStrip}.
$$
If $u(z_0)=a/r\Im z_0$, for some $z_0\in\MyDomain$, then $v(z_0)=0$ and hence, by the maximum principle for the subharmonic functions, $v(z)=0$ for each $z\in\MyDomain$. It follows that
$$
u(z)=\frac{\UpperFunc}{\UpperStrip}\Im z,
$$
and consequently, for each $x\in\R$,
$$
\limsup_{y\to0^+}\frac{u(x+iy)}{y}\leq\frac{\UpperFunc}{\UpperStrip}.
$$
Otherwise $v(z)<0$ for every $z\in\MyDomain$, 
that is
$$
u(z)<\frac{\UpperFunc}{\UpperStrip}\Im z,
$$
and for each $x\in\R$,
in view of Theorem (\ref{thm::SHHopfLemmaA}), 
$$
	\limsup_{y\to\UpperStrip^-}\frac{v(x+iy)}{y}<0,
$$
namely
$$
	\limsup_{y\to0^+}\frac{u(x+iy)}{y}<\frac{\UpperFunc}{\UpperStrip}.
$$
This proves the \statementSchwarz.
\qed
\section{\label{section:LaMetrica}Pseudo-metrics}
Let  $\CManA$ a (connected) complex manifold of dimension $n$.

Given $p\in\CManA$, $\TanVA\in T_p\CManA$,
we denote by $\Gamma_\CManA(p,\TanVA)$ the space of $C^1$ maps
$\gamma:]-\e,\e[\to\CManA$, for some $\e>0$, which satisfy $\gamma(0)=p$, 
$$
\gamma'(0):=d\gamma(0)(\frac{\rm d}{{\rm d}t})=\TanVA.
$$
For any subset $\SManA\sset\CManA$ we denote $\PSHU[\SManA,\CManA]$ the class of functions
$$
u:\CManA\to[0,\pi/4[
$$
which are plurisubharmonic in $\CManA$ and vanishing on $\SManA$.

As explained in the introduction, an element $\bar{u}\in\PSHU[\SManA,\CManA]$
is said to be {\it maximal} if $u(p)\leq\bar{u}(p)$
for every $u\in\PSHU[\SManA,\CManA]$ and $p\in\CManA$.
Clearly, a maximal element in $\PSHU[\SManA,\CManA]$ is unique,
provided it exists, and a maximal element exists in $\PSHU[\SManA,\CManA]$ if and only if
$$
\sup_{u\in\PSHU[\SManA,\CManA]}u\in\PSHU[\SManA,\CManA].
$$
Assume now that $\APair[\SManA,\CManA]$ is an analytic pair.
For $p\in\SManA$, $\TanVA\in T_p\SManA\sset T_p\CManA$, $u\in\PSHU[\SManA,\CManA]$ we set
$$
E_{u,\CManA}(p,\TanVA)=\inf_{\gamma\in\Gamma(p,\TanVA))}\limsup_{t\to0^+}\frac{u\bigl(\gamma(t)\bigr)}{t}
$$
and 
$$
E_{\SManA,\CManA}(p,\TanVA)=\sup_{u\in\PSHU[\SManA,\CManA]}E_{u,\CManA}(p,\TanVA).
$$
Clearly, if $\in\PSHU[\SManA,\CManA]$ admits a maximal element $\bar{u}$, one has
$$
E_{\bar{u},\CManA}(p,\TanVA)=E_{\SManA,\CManA}(p,\TanVA).
$$
Moreover, $E_{u,\CManA}$ and $E_{\SManA,\CManA}$ are positively homogeneous functions on $T\SManA$, i.e.
$$
E_{u,\CManA}(p,t\TanVA)=tE_{u,\CManA}(p,\TanVA),
$$
$$
E_{\SManA,\CManA}(p,t\TanVA)=tE_{\SManA,\CManA}(p,\TanVA).
$$
for $t>0$.
Observe that, in general, $E_{u,\CManA}$, $E_{\SManA,\CManA}$ are not symmetric with respect to
$\TanVA$.

Assuming a few of regularity on $u$ the definition of $E_{u,\CManA}(p,\TanVA)$ simplifies:
\begin{lemma}\label{lemma::ENoInfimum}
Let $u\in\PSHU[\SManA,\CManA]$ and $p\in\SManA$. If $u$ is Lipschitz in a neighbourhood of $p$ then for every $\TanVA\in T_p\SManA$ and $\gamma\in\Gamma_\CManA(p,J\TanVA)$
$$
E_{u,\CManA}(p,\TanVA)=\limsup_{t\to0^+}\frac{u\bigl(\gamma(t)\bigr)}{t}.
$$
\end{lemma}
\proof
It is sufficient to prove that for arbitrary $\gamma_1,\gamma_2\in\Gamma_\CManA(p,J\TanVA)$ it results
$$
\limsup_{t\to0^+}\frac{u\bigl(\gamma_1(t)\bigr)}{t}=\limsup_{t\to0^+}\frac{u\bigl(\gamma_2(t)\bigr)}{t}.
$$
Let $\gamma_1,\gamma_2\in\Gamma_\CManA(p,J\TanVA)$ and $z_1,\dots,z_n$ local complex coordinates near $x$. Then,
for $t>0$ sufficiently small, we have
\begin{eqnarray}
	u\bigl(\gamma_1(t)\bigr)
	&\leq&
	u\bigl(\gamma_2(t)\bigr)+\bigl(
		u\bigl(\gamma_1(t)\bigr)-u\bigl(\gamma_2(t)\bigr)\bigr)
	\nonumber\\
	&\leq&
	u\bigl(\gamma_2(t)\bigr)+C\bigl|z\bigl(\gamma_1(t)\bigr)-z\bigl(\gamma_2(t)\bigr)\bigr|
	\leq
	u\bigl(\gamma_2(t)\bigr)+o(t),
	\nonumber
\end{eqnarray}
and consequently
$$
\limsup_{t\to0^+}\frac{u\bigl(\gamma_1(t)\bigr)}{t}\leq\limsup_{t\to0^+}\frac{u\bigl(\gamma_2(t)\bigr)}{t}.
$$
Interchanging $\gamma_1$ and $\gamma_2$ we get the opposite inequality.
\qed
\begin{theorem}\label{TEO2}
Let 
$$
\CManA=\big\{z\in\C\mid \abs{\Im z}<\pi/4\big\}.
$$
Then the function $u(z)=\abs{\Im z}$ belongs to $\PSHU[\R,\CManA]$ and is maximal. Moreover, 
$$
E_{\R,\CManA}(x,\TanVA)=\abs{\TanVA}
$$
for every $x\in\R$, $\TanVA\in\R=T_x\R$.
\end{theorem}
\proof
Maximality is a consequence of Theorem (\ref{thm::PSHSchwarz}). Since $u$ is Lipschitz, Lemma (\ref{lemma::ENoInfimum}) then implies
$$
E_{\R,\CManA}(x,\TanVA)=\limsup_{t\to{0^+}}\frac{\abs{\Im{(x+it\TanVA)}}}{t}=\abs{\TanVA}.
$$
\qed
\begin{theorem}\label{thm:EDeltaC}
Let $\Delta$ be the (open) unit disc in $\C$, $I$ the interval $]-1,1[$. The function 
$$
u(z)=\abs{\Im(\arctanh(z))}
$$
 is maximal in $\PSHU[I,\Delta]$. Moreover, for every $x\in I$and $\TanVA\in\R=T_x I$ one has
$$
E_{I,\Delta}(x,\TanVA)=\frac{\abs{\TanVA}}{1-x^2}.
$$
\end{theorem}
\proof
We observe that the function
$$
f(z)=\arctanh(z)=\frac{1}{2}\log\frac{1+z}{1-z}
$$
is a biholomorphism between $\Delta$ and 
$$
\CManA=\bigl\{z\in\C\mid\abs{\Im(z)}<\pi/4\bigr\}
$$
 $f(I)=\R$ and 
$$
f'(z)=\frac{1}{1-z^2}.
$$
The statement is then an immediate consequence of Theorem \ref{thm::PSHSchwarz}.
\qed
The quantities $E_{\SManA,\CManA}$ decrease by holomorphic maps:
\begin{theorem}\label{thm::EMapdecrease}
Let $\APair[\SManA,\CManA]$, $\APair[\SManB,\CManB]$ be analytic pairs
and $f:\CManA\to\CManB$ a holomorphic map such that $f(\SManA)\sset\SManB$.
Then
\begin{equation}\nonumber
E_{\SManB,\CManB}\bigl(f(p),df(p)(\TanVA)\bigr)\leq E_{\SManA,\CManA}(p,\TanVA).
\end{equation}
for every $p\in\SManA$, $\TanVA\in T_p\SManA$,
\end{theorem}
\proof
We may assume that $E_{\SManA,\CManA}(p,\TanVA)<+\infty$. Let $\e>0$ be fixed and $u\in\PSHU[\SManB,\CManB]$. Then $u\circ f\in\PSHU[\SManA,\CManA]$ and by definition of $E_{u,\CManA}$ there exists $\gamma\in\Gamma_\CManA(p,J\TanVA)$ such that
$$
\limsup_{t\to0^+}\frac{u\circ f\bigl(\gamma(t)\bigr)}{t}<E_{u\circ f,\CManA}(p,\xi)+\e\leq E_{\SManA,\CManA}(p,\TanVA)+\e.
$$
Then $u\circ f\in\Gamma_\CManB\bigl(f(p),Jdf(p)(\TanVA)\bigr)$ and consequently
\begin{eqnarray}
E_{u,\CManB}\bigl(f(p),df(p)(\TanVA)\bigr)&\leq&\limsup_{t\to0^+}\frac{u\bigl(f\circ\gamma(t)\bigr)}{t}\nonumber\\
&=&\limsup_{t\to0^+}\frac{u\circ f\bigl(\gamma(t)\bigr)}{t}<E_{\SManA,\CManA}(p,\TanVA)+\e.\nonumber
\end{eqnarray}
Since $\e>0$ is arbitrary we get 
$$
E_{u,\CManB}\bigl(f(p),df(p)(\TanVA)\bigr)\leq E_{\SManA,\CManA}(p,\TanVA).
$$
We obtain the desired inequality taking the supremum over $\PSHU[\SManB,\CManB]$
\qed
Consider now the unit disc $\Delta$ and recall that for $p\in\SManA$, $\TanVA\in T_p\SManA$ one defines
\begin{eqnarray}
F_{\SManA,\CManA}(p,\TanVA)=\inf\bigl\{a>0&\mid&\exists f\in\Hol(\Delta,\CManA),f(]-1,1[\sset\CManA,\nonumber\\
&&f(0)=p,f'(0)=a^{-1}\TanVA\bigr\}.
\end{eqnarray}
(cf. \cite{article:GTV})
If $f\in\Hol(\Delta,\CManA)$ satisfies $f(]-1,1[\sset\CManA,f(0)=p,f'(0)=a^{-1}\TanVA$ then, in view of Theorems \ref{thm::EMapdecrease}, \ref{thm:EDeltaC}, we have
$$
E_{\SManA,\CManA}(p,\TanVA)\leq
E_{]-1,1[,\Delta}(0,a)=a;
$$
taking the infimum of $a$ over all maps $f\in\Hol(\Delta,\CManA)$ we get:
\begin{theorem}\label{thm::ELessThanF}
Let $\SManA\sset\CManA$ be an analytic pair. Then
$$
E_{\SManA,\CManA}(p,\TanVA)\leq F_{\SManA,\CManA}(p,\TanVA).
$$
for every $p\in\SManA$, $\TanVA\in T_p\SManA.$
\end{theorem}
The theorem which follows characterizes the ``complex geodesic''
for the pseudo-metrics $E_{\SManA,\CManA}$.
\begin{theorem}\label{thm::EGeodesics}
Let $\APair[\SManA,\CManA]$ be an analytic pair.
Let $S=\bigl\{\abs{\Im(z)}<\pi/4\bigr\}$ and $f:S\to\CManA$ be a
holomorphic map such that $f(\R)\sset\SManA$.
Then for a function $u\in\PSHU[\SManA,\CManA]$ the following conditions are equivalent:
\begin{itemize}
\item[i)] for every $z\in S$
$$
u\bigl(f(z)\bigr)=\abs{\Im(z)};
$$
\item[ii)] for every $x\in\R$, $\TanVA\in\R=T_x\R$
$$
E_{u,\CManA}\bigl(f(x),df(x)(\TanVA)\bigr)=\abs{\TanVA};
$$
\item[iii)] there is $x_0\in\R$ such that 
$$
E_{u,\CManA}\bigl(f(x_0),df(x_0)(\TanVA)\bigr)=\abs{\TanVA}.
$$
for every $\TanVA\in\R=T_x\R$.
\end{itemize}
Moreover, if such conditions are fulfilled, for every $x\in\R$, $\TanVA\in\R$ the following identities hold
$$
	E_{u,\CManA}\bigl(f(x),df(x)(\TanVA)\bigr)
	=E_{\SManA,\CManA}\bigl(f(x),df(x)(\TanVA)\bigr)
	=F_{\SManA,\CManA}\bigl(f(x),df(x)(\TanVA)\bigr)
	=\abs{\TanVA}.
$$
\end{theorem}
\proof
The implications
i)$\implies$ ii), ii)$\implies$ iii) are evident and that iii)$\implies$i) follows immediately from Theorem \ref{thm::PSHSchwarz}. In order to prove the last equality it is sufficient to observe that by definition 
$$
E_{u,\CManA}\bigl(f(x),df(x)(\TanVA)\bigr)\leq E_{\SManA,\CManA}\bigl(f(x),df(x)(\TanVA)\bigr);
$$
moreover, by Theorem \ref{thm::ELessThanF}
$$
E_{\SManA,\CManA}\bigl(f(x),df(x)(\TanVA)\bigr)\leq F_{\SManA,\CManA}\bigl(f(x),df(x)(\TanVA)\bigr),
$$
and by the properties of $F_{\SManA,\CManA}$ (cf. \cite{article:GTV})
$$
F_{\SManA,\CManA}\bigl(f(x),df(x)(\TanVA)\bigr)\leq F_{\R,S}(x,\TanVA)=\abs{\TanVA}.
$$
Then if ii) holds 
$$
\abs{\TanVA}=E_{u,\CManA}\bigl(f(x),df(x)(\TanVA)\bigr).
$$
\qed
The holomorphic maps $f:S\to\CManA$ which satisfy $f(\R)\sset\SManA$ and the conditions of Theorem \ref{thm::EGeodesics} are called {\sl $E_{u}-$complex geodesic}.
The $E-$complex geodesic are a useful tool to give sufficient conditions in order to state maximality of plurisubharmonic functions in $\PSHU[\SManB,\CManB]$.
\begin{theorem}
Let $(\SManA,\CManA)$ an analytic pair and $u\in\PSHU[\SManA,\CManA]$. Suppose that for every $q\in\CManA\setminus\SManA$ there esists an $E_{u}-$complex geodesic $f:S\to\CManA$ such that $q\in f(S)$. Then $u$ is maximal.
\end{theorem}
\proof
Let $w\in\PSHU[\SManA,\CManA]$. We have to prove that $w(q)\leq u(q)$ for every $q\in\CManA$ so let $q\in\CManA$. If $q\in\SManA$ then
$w(q)=u(q)=0$ and in such a case the thesis is evident, so we assume that $q\in\CManA\setminus\SManA$. Let $S=\bigl\{\abs{\Im(z)}<\pi/4\bigr\}$ and $f:S\to\CManA$ una $E-$complex geodesic such that $f(z_0)=q$, $z_0\in S$. In view of Theorem\ref{thm::EGeodesics} we have $u(q)=\abs{\Im(z_0)}$. Observe now that $u\bigl(f(z)\bigr)$ is subharmonic in $S$ and $0\le u\bigl(f(z)\bigr)\le \pi/4$, so, in view of Theorem \ref{thm::PSHSchwarz} we have $w\bigl(f(z)\bigr)\leq\abs{\Im(z)}$ for every $z\in S$. In particular
$$
v(q)=v\bigl(f(z_0)\bigr)\leq\abs{\Im(z_0)}=u(q).
$$
and from this it follows that $u$ is maximal, $q$ being arbitrary.
\qed
\section{\label{section:MongeAmpere}Complex Monge-\Ampere\ equation}
The theorem which follows put in evidence  the relationship between maximal functions for analytic pairs $(V,M)$ and solutions of the complex Monge-\Ampere\ equation on $M\setminus V$.
For the main results about existence, unicity and maximum principle for solutions of the complex Monge-\Ampere\ equation we refer to \cite{article:BedfordMASurvey}.
\begin{theorem}\label{thm::MAMaximal}
Let $(\SManA,\CManA)$ be an analytic pair and $u\in\PSHU[\SManA,\CManA]$. Then
\begin{itemize}
\item[i)] if $u$ is continuous and maximal, $(dd^cu)^n=0$ on 
$\CManA\setminus\SManA$;\\
\item[ii)] if $u$ is an exhaustion function such that $(dd^cu)^n=0$ on 
$\CManA\setminus\SManA$ and $u(p)=0$ if and only $p\in\SManA$ then $u$ is maximal.
\end{itemize}
In particular, if $u$ is a continuous exhaustion function and $u(p)=0$ if and only if $p\in\SManA$
then $u$ is maximal if and only if $(dd^cu)^n=0$ on $\CManA\setminus\SManA$.
\end{theorem}
\proof
Assume that $u$ is continuous and maximal and let us show that actually $u$ is a solution of $(dd^cu)^n=0$ on $\CManA\setminus\SManA$.

Let $p\in\CManA\setminus\SManA$ and $U\sset\CManA\setminus\SManA$
be a relatively compact neighbourhood of $p$.
Let $w:\CManA\to[0,\pi/4]$ be the function defined by
$$
w(p)=\left\{\begin{array}{ll}u(p)\quad&p\in\CManA\setminus U\\ 
v(p)\quad&p\in\bar{U}
\end{array}\right.,
$$
where $v$ is the solution of the problem
$$
\left\{\begin{array}{ll}(dd^cv)^n=0\quad&{\rm in\ } U\\
v=u\quad&{\rm in\ } \partial U.
\end{array}
\right.
$$
The function $v$ is characterized by
$$
v(p)=\sup\{w(p)\}
$$
where the supremum is taken over the set of functions $w$ which are plurisubharmonic in $U$, continuous on $\bar{U}$ and satisfying $w\leq u$ on $\partial U$. The function $w$ belongs to $\PSHU[\SManA,\CManA]$ and, by construction, $w\ge u$. Since $u$ is maximal then  $u=w$; in particular $u$ is a solution of Monge-\Ampere\ in a neighbourhood of $p\in\CManA\setminus\SManA$. Thus $u$ is a solution of $(dd^cu)^n=0$
on $\SManA\setminus\CManA$, $p\in\CManA\setminus\SManA$ being arbitrary.

Conversely, suppose that $u$ is an exhaustion function for $\CManA$ and a solution of the Monge-\Ampere\  equation
on $\SManA\setminus\CManA$. In particular, $\SManA$ is a compact submanifold of $\CManA$.
Let $w$ be an arbitrary function of $\in\PSHU[\SManA,\CManA]$. We have to prove that $w(p)\le u(p)$ for every $p\in M$. This is certainly true if $p\in\CManA$ since then 
$w(p)=u(p)=0$. Thus we assume that $p\in\CManA\setminus\SManA$. By hypothesis $u(p)>0$. Let $\e>0$ be such that $u(p)<\frac{\pi}{4}-\e$,
$$
D=\bigl\{q\in\SManA\mid0<u(q)<\frac{\pi}{4}-\e\bigr\};
$$
since $u$ is an exhaustion function $D$ is relatively compact. Let $F_1$ be the subset of the boundary $\partial D$ of $D$ where $u$ takes the value $\frac{\pi}{4}-\e$
$$
u_{\e}=\frac{\pi/4}{\pi/4-\e}u.
$$
Then $\partial D=M\cup F_1$ and we are going to show that $u_{\e}\ge w$ on $\partial D=M\cup F_1$. Indeed, if $q\in\SManA$ then $w(q)=0=u_\e(q)$ and if $q\in F_1$ then $w(q)<\pi/4=u_\e(q)$. Since $w$ is plurisubharmonic we have $(dd^cv)^n\geq0=(dd^cu_\e)^n$
on $D$, whence $w(q)\leq u_\e(q)$ for every $q\in D$. In particular $w(p)\leq u_\e(p)$ for every $\e>0$, hence $w(p)\leq u(p)$. Since $p\in\CManA$ is arbitrary $u\ge w$ on $\CManA$.
\qed
\section{\label{section:CasoLiscio}The smooth case}
Let $\CManA$ be a complex manifold and $u$ a $C^2$ function on $\CManA$.
Let $p\in\CManA$, $\TanVA\in T_p\CManA$ and $f$ be a germ at $0\in\C$ of a holomorphic map with values in $\CManA$  such that $f(0)=p$, $f'(0)=\TanVA$. Then the complex number
$$
\Levi_u(p,\TanVA)=\frac{\partial^2(u\circ f)(0)}{\partial z\partial\bar{z}},
$$
depends only on $u, p, \TanVA$ and it is nothing but that the Levi form of $u$ at $p$ evaluated at $\TanVA$. 
\begin{proposition}\label{EG}
Let $\APair[\SManA,\CManA]$ be an analytic pair and $u\in\PSHU[\SManA,\CManA]$. Assume that $u^2$ is $C^2$ around $V$. Then, for every $p\in\SManA$ and $\TanVA\in T_p\SManA$, we have the following equality
$$
E_{u,\CManA}(p,\TanVA)=\sqrt{2\Levi_{u^2}(p,\TanVA)}.
$$
\end{proposition}
\proof
Let $p\in\SManA$, $\TanVA\in T_p\SManA$ and $f$ be a holomorphic map with values in $\CManA$, defined in a neighbourhood $U$ of the origin $0\in\C$ and such that $f(0)=p$, $f'(0)=\TanVA$. Since $u^2$ is of class $C^2$ in a neighbourhood of $\SManA$ it is locally Lipschitz in a neighbourhood of $\SManA$. Then Lemma \ref{lemma::ENoInfimum} implies
$$
E_{u,\CManA}(p,\TanVA)=\limsup_{y\to0^+}\frac{u\circ f(iy)}{\TanVA}.
$$
Now we set
$$
g(y)=\bigl(u\circ f(iy)\bigr)^2,
$$
and observe that, since $u\circ f$ vanishes on
$U\cap\R$, one has
$$
2\Levi_{u^2}(p,\TanVA)=
2\frac{\partial^2(u\circ f)^2}{\partial z\partial\bar{z}}(0)=
\frac{1}{2}\Delta(u\circ f)^2(0)=
\frac{1}{2}g''(0),
$$
where
$$
\Delta=\frac{\partial^2}{\partial x^2}+\frac{\partial^2}{\partial y^2}.
$$
On the other hand, since $g(y)$ is $C^2$, non negative and vanishing on $y=0$, from the elementary identity
$$
\sqrt{\frac{1}{2}g''(0)}=
\lim_{y\to0}\frac{\sqrt{g(y)}}{y}
$$
we get 
$$
\sqrt{2\Levi_{u^2}(p,v)}=\sqrt{\frac{1}{2}g''(0)}
=\lim_{y\to0^+}\frac{u\circ f(iy)}{y}
=E_{u,\CManA}(p,\TanVA).
$$
This proves the proposition.
\qed
\begin{theorem}\label{LI}
Let $\APair[\SManA,\CManA]$ be an analytic pair and $u\in\PSHU[\SManA,\CManA]$.
Assume that $(\SManA,\CManA,u)$ is a Monge-Amp\`ere model. Then, for every $p\in\SManA$, $\TanVA\in T_p\SManA$
$$
E_{\SManA,\CManA}(p,\TanVA)=F_{\SManA,\CManA}(p,\TanVA)=\sqrt{2\Levi_{u^2}(p,\TanVA)}.
$$
In particular, $(\SManA,\CManA,u)$ is a maximal model.
\end{theorem}
\proof
We have that $\SManA=\{u=0\}$, $(dd^cu)^n=0$ on $\CManA\setminus\SManA$ and $u^2$ is a smooth, strictly plurisubharmonic exhaustion function for $\CManA$.

Let $p\in\SManA$, $\TanVA\in T_p\SManA$. Without loss of generality we may assume 
$$
\sqrt{2\Levi_{u^2}(p,\TanVA)}=1.
$$ 
In view of Theorem \ref{thm::MAMaximal} the function $u\in\PSHU[\SManA,\CManA]$ is maximal i.e.
$$
E_{\SManA,\CManA}(p,\TanVA)=E_{u,\CManA}(p,\TanVA),
$$ 
so, by Proposition \ref{EG},
$$
E_{u,\CManA}(p,\TanVA)=\sqrt{2\Levi_{u^2}(p,\TanVA)}=1.
$$
Let us denote $g$ the Riemannian metric induced on $V$ by the restriction to $TV$ of the Levi form $dd^c(u^2)$ 
(cf. \cite{article:PatrizioWong}). Since $u$ is an exhaustion function for $M$, $(V,g)$ is a compact Riemannian manifold, so there is a geodesic $\g:\R\to\CManA$ such that $\g(0)=p$, $\g'(0)=\TanVA$ and
$$
g\bigl(\g(x),\g'(x),\g'(x)\bigr)=1
$$
for every $x\in\R$.

In view of the results proved in \cite{article:PatrizioWong}
(cf. also Theorem 3.1 of \cite{article:LempertSzoke}),
if $S=\bigl\{\abs{\Im z}<\pi/4\bigr\}$, the map $f:S\to\CManA$ defined by 
$$
f(z)=f(x+iy)=d\g(x)(y)
$$
is holomorphic.
By construction $f(\R)\sset\CManA$ and, moreover,
$$
1=E_{u,\CManA}(p,\TanVA)=E_{u,\CManA}(f(0),f'(0)).
$$
Theorem \ref{thm::EGeodesics} now implies that $f$
is an $E-$complex geodesic and, since $p=f(0)$, $\TanVA=f'(0)$, one has
$$
F_{\SManA,\CManA}(p,\TanVA)=
E_{\SManA,\CManA}(p,\TanVA)=
E_{u,\CManA}(p,\TanVA)=\sqrt{2\Levi_{u^2}(p,\TanVA)}=1.
$$
\qed
\section{\label{section:esa}Convex homogeneous real functions}
For every point $z\in\C^n$ we set $z=x+iy$, $x,y\in\R^n$, $x=\Re z$, $y=\Im z$.\\
Let $\mink:\R^n\to[0,+\infty[$ be a positively homogeneous convex function such that $\mink(x)=0$ if and only if $x=0$; $\mink$ is the Minkowsky functional associated to a bounded open convex subset of $\R^n$.
Observe that we do not require the property $\mink(-x)=\mink(x)$ for $x\in\R^n$.

Let 
$$
X_\mink=\bigl\{z\in\C^n\mid \mink\bigl(\Im(z)\bigr)<\pi/4\bigl\}
$$
and $u_\mink:X_\mink\to[0,\pi/4[$ be the function defined by
$$
u_\mink(z)=\mink\bigl(\Im(z)\bigr)
$$
for every $z\in X_\mink$.
\begin{theorem}
Let $X_\mink$, $u_\mink$ be as above. Then $u_\mink\in\PSHU[\R^n,X_\mink]$ is maximal and for every $x\in\R^n$, $\xi\in\R^n=T_x\R^n$ the following identity holds true
$$
E_{\R^n,X_\mink}(x,\xi)=E_{u_\mink,X_\mink}(x,\xi)=\mink(\xi).
$$
\end{theorem}
\proof
Assume first that $\mink$ is of class $C^2$. Then, since the function
$u_\mink(x+iy)$ does not depend on $x$ it follows that
$$
\frac{\partial^2u_\mink(x+iy)}{\partial z_i\partial\bar{z}_j}=
\frac{1}{4}\frac{\partial^2\mink(y)}{\partial y_i\partial y_j};
$$
it follows that $u_\mink$ is plurisubharmonic, since $\mink$ is convex.
If $\mink$ is only continuous the same conclusion is obtained approximating $\mink$ by smooth functions.

In order to show the maximality of $u_\mink$, we have to prove that, if $w\in\PSHU[\R^n,X_\mink]$, then $w(z)\leq u_\mink(z)$ for every $z\in X_\mink$.

This is obviously true if $\Im z=0$ for then $w(z)=0\leq u_\mink(z)$, so let $\Im z\neq 0$. Define on 
$$
S^+=\bigl\{\z\in\C\mid0<\Im\z<\pi/4\big\}
$$
the function $f:S^+\to X_\mink$ setting, for every $\z\in\C$
$$
f(\z)=\Re z+\z\mink(\Im z)^{-1}\Im z.
$$
Then, by construction
$$
u_\mink\bigl(f(\z)\bigr)=\mink\bigl((\Im\z)\mink(\Im z)^{-1}\Im z\bigr)=\Im\z.
$$
The function $w\circ f:S^+\to[0,\pi/4[$ is subharmonic and satisfies
$$
\limsup_{S\ni z\to x}w\circ f(z)=0
$$
for every $x\in\R$, so, in view of Theorem \ref{thm::PSHSchwarz}, we have
$$
w\bigl(f(\z)\bigr)\leq\Im\z=u_\mink\bigl(f(\z)\bigr).
$$
for every $\z\in S^+$. In particular, for $\z_0=i\mink(\Im z)$ we get $f(\z_0)=z$ and consequently
$$
v(z)=v\bigl(f(\z_0))\bigr)\leq u_\mink\bigl(f(\z_0))\bigr)=u_\mink(z).
$$
This proves that $u_\mink$ is maximal.

Now we observe that, since $\mink$ is a convex function, $u_\mink$ is Lipschitz; then, by Lemma \ref{lemma::ENoInfimum}, we have
$$
E_{\R^n,X_\mink}(x,\TanVA)=E_{u_\mink,X_\mink}(x,\TanVA)=
\lim_{t\to0^+}\frac{u_\mink(x+ity)}{t}=
\lim_{t\to0^+}\frac{\mink(ty)}{t}=\mink(y),
$$
for every $x\in\R^n$, $\TanVA\in\R^n=T_x\R^n$.\\ 
This proves completely the theorem.
\qed
\begin{remark}\label{MODELLI}
{\rm
It would be interesting to provide a characterization of the 
models $(\R^n, X_\mink, u_\mink)$ as done by Abate e Patrizio in \cite{article:AbatePatrizioConstCurv}, where $\mink$ is assumed to be in $C^{\infty}(\R^n\setminus\{0\})$ and symmetric i.e. $\mink(-x)=\mink(x)$ for every $x\in\R^n$.}
\end{remark}
\section{\label{section:TuboEllittico}The elliptic tube of Lempert}
{
The following construction is due to Lempert \cite{article:LempertEllipticTube}.

Given a segment $I\sset\R^n\sset\C^n$ of positive length we denote $L(I)\sset\C^n$ the unique complex line which contains $I$. We assume that $I$ is a relatively open interval in the real straight line of $L(I)$ containing it. Let  $\tilde{I}\sset L(I)$ be the relatively open disc in $L(I)$ whose diameter is $I$.

Let now $\Cvx\sset\R^n$ be a convex domain. The {\sl elliptic tube} over $\Cvx$ is defined by
$$
\Cvx^{\rm ell}=\bigcup\limits_ I\bigl\{\tilde{I}\mid I\sset \Cvx, I\ {\rm segment}\bigr\}.
$$
The main result of this section consists of finding the maximal function $u\in\PSHU[\Cvx,\Cvx^{\rm ell}]$ and the explicit computation of $E_{\Cvx,\Cvx^{\rm ell}}$ when $D$ is a bounded convex domain in $\R^n$ .

If $z=x+iy\in\C^n$, $x\in\Cvx$, consider the functional of Minkowski of $D$ centered at $x$ and evaluated at $y$
$$
p(z)=p_\Cvx(z)=\inf\bigl\{t>0\mid x+t^{-1}y\in\Cvx\bigr\}
$$
It is easy to check that a point $z=x+iy\in\C^n$ belongs to $z\in\Cvx^{\rm ell}$ if and only if $x=\Re z\in\Cvx$ and $p(z)p(\bar{z})<1$.

We have the following 
\begin{theorem}\label{TUBO}
Let $\Cvx\sset\R^n$ be a bounded convex domain and $p(z)$ the Minkowski functional. Let $u:\Cvx^{\rm ell}\to[0,\pi/4[$ be defined by
\begin{equation}\label{eq::TubeMAFunction}
u(z)=u_\Cvx(z)=\frac{\arctan\bigl(p(z)\bigr) + \arctan\bigl(p(\bar z)\bigr))}{2}.
\end{equation}
Then
\begin{itemize}
\item[1)] $u$ is locally Lipschitz and $u(z)=u(x+iy)=0$ if and only if $y=0$;
\item[2)] $u$ is plurisubharmonic on $\Cvx^{\rm ell}$;
\item[3)] $u\in\PSHU[\Cvx,\Cvx^{\rm ell}]$ is maximal;
\item[4)] $(dd^c u)^n=0$ on $\Cvx^{\rm ell}\setminus\Cvx$.
\end{itemize}
Moreover, if $z=x+iy$ with $x\in\Cvx$, $y\in\R^n=T_x\Cvx$, then
$$
E_{u,\Cvx^\e}(x,y)=
E_{\Cvx,\Cvx^\e}(x,y)=
F_{\Cvx,\Cvx^\e}(x,y)=\frac{p(z)+p(\bar{z})}{2}.
$$
Finally, if $\partial\Cvx$ is of class $C^2$ then also $p$ and $u$ are of class $C^2$
on $\Cvx^{\rm ell}\setminus\Cvx$ and for each $i,j=1,\ldots,n$
\begin{equation}\label{eq::TubeLevi}
\frac{\partial^2u(z)}{\partial z_i\partial\bar{z}_j}=
\frac{1}{4}\left(\frac{\partial^2p(z)}{\partial y_i\partial y_j}
+\frac{\partial^2p(\bar{z})}{\partial y_i\partial y_j}\right).
\end{equation}
\end{theorem}
\proof
The statement 1) is a consequence of convexity and boundedness of $\Cvx$.

In order to prove 2) and 3) define for $z\in\Cvx$ 
$$
\tilde{u}(z)=\sup\bigl\{v(z)\mid v\in\PSHU[\Cvx, \Cvx^\e]\bigr\};
$$
it is then sufficient to show that   
\begin{equation}\nonumber
u(z)=\tilde{u}(z)
\end{equation}
for every $z\in\Cvx^\e$.
If  $z\in\Cvx$ the equality is evident, so we assume that $z=x+iy\in\Cvx^\e$ with $y\neq0$.
Set
\begin{eqnarray*}
&&t_1=p(z)^{-1}\\
&& t_2=p(\bar{z})^{-1}\\
&&x_1=x+t_1y\\
&& x_2=x-t_2y.
\end{eqnarray*}
Then $x_1,x_2\in\partial\Cvx$ and the segment with endpoints
$x_1$ and $x_2$ is contained in $\Cvx$; we easily derive
\begin{eqnarray}
&&x=\frac{t_2}{t_1+t_2}x_1+\frac{t_1}{t_1+t_2}x_2,\nonumber\\&&y=\frac{1}{t_1+t_2}x_1-\frac{1}{t_1+t_2}x_2,\nonumber
\end{eqnarray}
namely
$$
z=\frac{t_2+i}{t_1+t_2}x_1+\frac{t_1-i}{t_1+t_2}x_2.
$$

Let $\Delta$ be the unit disc in $\C$ and $f:\Delta\to\C^n$  defined by
$\z\in\Delta$
$$
f(\z)=\frac{1-\z}{2}x_1+\frac{1+\z}{2}x_2.
$$
Since $f$ sends $]-1,1[$ in $D$ then $f(\Delta)\sset\Cvx^{\rm ell}$ (\cite{article:LempertEllipticTube}). Setting
$$
\z_0=\frac{t_1-t_2}{t_1+t_2}+i\frac{2}{t_1+t_2};
$$
the inequality $p(z)p(\bar{z})<1$ implies $\z_0\in\Delta$; moreover $f(\z_0)=z$ so
\begin{eqnarray}
\arctanh(\z_0)&=&\frac{1}{2}\log\frac{1+\z_0}{1-\z_0}=\frac{1}{2}\log\frac{t_1-i}{t_2+i}\nonumber\\
&=&\frac{1}{2}\log\left(\frac{t_1t_2-1}{t_2^2+1}-i\frac{t_1+t_2}{t_2^2+1}\right)\nonumber
\end{eqnarray}
whence
\begin{eqnarray*}
\left|\Im\bigl(\arctanh(\z_0)\bigr)\right|\!\!\!
&=&\!\!\!\frac{1}{2}\left|\arg\left(\frac{t_1t_2-1}{t_2^2+1}-i\frac{t_1+t_2}{t_2^2+1}\right)\right|=\\
&&\!\!\!\frac{1}{2}\arctan\frac{t_1+t_2}{t_1t_2-1}.
\end{eqnarray*}
Since $t_1=p(z)^{-1}$, $t_2=p(\bar z)^{-1}$ we deduce
$$
\left|\Im\bigl(\arctanh(\z_0)\bigr)\right|=
\frac{1}{2}\arctan\frac{p(z)+p(\bar z)}{1-p(z)p(\bar z)}.
$$
Finally, putting
\begin{eqnarray*}
&&p(z)=\tan(\arctan(p(z))\\
&& p(\bar z)=\tan(\arctan(p(\bar z))
\end{eqnarray*}
into the elementary formula
$$
\tan(\a+\bg)=\frac{\tan(\a)+\tan(\bg)}{1-\tan(\a)\tan(\bg)},
$$
we get
$$
\left|\Im\bigl(\arctanh(\z_0)\bigr)\right|=\frac{\arctan\bigl(p(z)\bigr)+\arctan\bigl(p(\bar z)\bigr)}{2}=u(z).
$$
Now let $w\in\PSHU[\Cvx,\Cvx^{\rm ell}]$. Then $w\circ f:\Delta\to[0,\pi/4[$ is subharmonic and vanishing on $]-1,1[$. In view of Theorem \ref{thm:EDeltaC}, for every $\z\in\Delta$ we have
$$
w\circ f(\z)\leq\left|\Im\bigl(\arctanh(\z)\bigr)\right|;
$$
in particular
$$
v(z)=v\circ f(\z_0)\leq\left|\Im\bigl(\arctanh(\z_0)\bigr)\right|=u(z)
$$
since $w$ is arbitrary
$$
\tilde{u}(z)\leq u(z).
$$
As for the opposite inequality we observe that there is a holomorphic map $g:\Cvx^{\rm ell}\to\Delta$ saisfying $g(\Cvx)\sset]-1,1[$, $g(f(\z))=\z$ for every $\z\in\Delta$ (cf. \cite{article:LempertEllipticTube}). It follows
$$
\left|\Im\bigl(\arctanh(g)\bigr)\right|\in\PSHU[\Cvx,\Cvx^\e],
$$
which implies
$$
\tilde{u}(z)\geq\left|\Im\bigl(\arctanh(g(z))\bigr)\right|
=\left|\Im\bigl(\arctanh(\z_0))\bigr)\right|
=u(z).
$$
since $g(z)=g\bigl(f(\z_0)\bigr)=\z_0$,

4) is an immediate consequence of continuity of $u$ in view of 
Theorem \ref{thm::MAMaximal}.

In order to prove the last equality let $x\in\Cvx$, $y\in\R^n=T_x\Cvx$. Since $u$ is locally Lipschitz, by Lemma \ref{lemma::ENoInfimum} we have
 \begin{eqnarray}
E_{u,\Cvx^\e}(x,y)=E_{\Cvx,\Cvx^\e}(x,y)&=&
\lim_{t\to0^+}\frac{\arctan\bigl(tp(z)\bigr)+\arctan\bigl(tp(\bar z)\bigr)}{2t}\nonumber\\
&=&\frac{p(z)+p(\bar{z})}{2},\nonumber
\end{eqnarray}
where $z=x+iy$,
Let $S=\bigr\{\abs{\Im(z)}<\pi/4\bigl\}$. We are going to show that the map $h:S\to D^{\rm ell}$ defined by
$$
h(\eta)=f\bigr(\tanh(\eta)\bigl).
$$
is an $E-$complex geodesic.

Define $k:\Cvx^{\rm ell}\to S$ by
$$
k(z)=\arctanh\bigl(g(z)\bigr).
$$
In view of the identity $g\circ f(\z)=\z$ we get $k\circ h(\eta)=\eta$ for every $\eta\in S$; it follows 
\begin{eqnarray}
\abs{\xi}\!\!\!&=&\!\!\!
E_{\R,S}(x,\xi)\leq
E_{\Cvx,\Cvx^\e}(h(x),dh(x)(\xi)\nonumber\\
&\leq&\!\!\!E_{\R,S}(k\circ h(x),d(k\circ h)(x)(\xi)=
E_{\R,S}(x,\xi)=\abs{\xi}\nonumber
\end{eqnarray}
for every $x\in\R$, $v\in\R=T_x\R$, and consequently
$$
E_{\Cvx,\Cvx^{\rm ell}}(h(x),dh(x)(\xi)=\abs{\xi}.
$$
Theorem \ref{thm::EGeodesics} now implies that $h$ is a $E-$complex geodesic. Moreover, again in view of Theorem \ref{thm::EGeodesics}, we have
\begin{eqnarray}
E_{\Cvx,\Cvx^{\rm ell}}(x,v)\!\!\!&=&\!\!\!E_{\Cvx,\Cvx^{\rm ell}}(h(k(x)),dh(k(x))(dk(x)(v)))\nonumber\\
&=&\!\!\!F_{\Cvx,\Cvx^{\rm ell}}(h(k(x)),dh(k(x))(dk(x)(v)))=F_{\Cvx,\Cvx^{\rm ell}}(x,v).\nonumber
\end{eqnarray}
\def\varxi{x_i}
\def\varxj{x_j}
\def\varyi{y_i}
\def\varyj{y_j}
\def\varzi{z_i}
\def\varbzj{\bar{z}_j}
\def\myP{p}
\def\auxF{\mu}
\def\DiniF{U}
\def\myIndexA{\alpha}
\def\myIndexB{\beta}
\def\mySum{\sum_{\myIndexA=1}^n}
\def\myDoubleSum{\sum_{\myIndexA,\myIndexB=1}^n}

\def\myDerA[#1,#2]{\frac{\partial{#1}}{\partial{#2}}}
\def\myDerAA[#1,#2,#3]{\frac{\partial^2{#1}}{\partial{#2}\partial{#3}}}
\def\myDerB[#1,#2]{{#1}_{#2}}
\def\myDerBB[#1,#2,#3]{{#1}_{{#2}{#3}}}
\def\myDerC[#1]{D_{#1}}

Assume now that $\partial\Cvx$ is smooth of class $C^2$.
Then we claim that the function $p(z)$ (and hence $u(z)$)
is smooth of class $C^2$ on $\Cvx^{\rm ell}\setminus\Cvx$
and for $i,j=1,\ldots,n$,
\begin{eqnarray}
	\label{eq::RealDX}
	&&\myDerA[\myP,\varxi]=\myP\myDerA[\myP,\varyi],
	\\
	\label{eq::RealDXY}
	&&\myDerAA[\myP,\varxi,\varyj]=
		\myP\myDerAA[\myP,\varyi,\varyj]+\myDerA[\myP,\varyi]\myDerA[\myP,\varyj],
	\\
	\label{eq::RealDXX}
	&&\myDerAA[\myP,\varxi,\varxj]=
		\myP^2\myDerAA[\myP,\varyi,\varyj]+2\myP\myDerA[\myP,\varyi]\myDerA[\myP,\varyj].
\end{eqnarray}
Assuming for granted such relations we have
\begin{eqnarray}\label{eq::DTwoArctan}
	\myDerAA[\arctan\myP,\varzi,\varbzj]
	=(1+\myP^2)^{-2}\left[
		(1+\myP^2)\myDerAA[\myP,\varzi,\varbzj]
		-2\myP\myDerA[\myP,\varzi]\myDerA[\myP,\varbzj]
	\right].
\end{eqnarray}
Since
\begin{eqnarray}
	\myDerA[\phantom{\myP},\varzi]
	&=&\frac{1}{2}\left(\myDerA[\phantom{\myP},\varxi]+\frac{1}{i}\myDerA[\phantom{\myP},\varyi]\right),
	\nonumber\\
	\myDerA[\phantom{\myP},\varbzj]
	&=&\frac{1}{2}\left(\myDerA[\phantom{\myP},\varxj]-\frac{1}{i}\myDerA[\phantom{\myP},\varyj]\right),
	\nonumber
\end{eqnarray}
then, using (\ref{eq::RealDX}) we obtain
\begin{eqnarray}
	\myDerA[\myP,\varzi]
	&=&\frac{1}{2}\left(\myP-i\right)\myDerA[\myP,\varyi],
	\nonumber\\
	\myDerA[\myP,\varbzj]
	&=&\frac{1}{2}\left(\myP+i\right)\myDerA[\myP,\varyj],
	\nonumber
\end{eqnarray}
and hence
\begin{eqnarray}\label{eq::RelOne}
	&&\myDerA[\myP,\varzi]\myDerA[\myP,\varbzj]
	=\frac{1}{4}\left(1+\myP^2\right)\myDerA[\myP,\varyi]\myDerA[\myP,\varyj].
\end{eqnarray}
Recalling that
\begin{eqnarray}
	\myDerAA[\phantom{\myP},\varzi,\varbzj]
	&=&\frac{1}{4}\left[
			\left(
				\myDerAA[\phantom{\myP},\varxi,\varxj]
				+\myDerAA[\phantom{\myP},\varyi,\varyj]
			\right)+
			\frac{1}{i}\left(
				\myDerAA[\phantom{\myP},\varxi,\varyj]
				-\myDerAA[\phantom{\myP},\varxj,\varyi]
			\right)
		\right],
	\nonumber\\
	\nonumber
\end{eqnarray}
from (\ref{eq::RealDXY}) and (\ref{eq::RealDXX}) we obtain
\begin{eqnarray}
	\myDerAA[\myP,\varzi,\varbzj]
	&=&\frac{1}{4}\left(1+\myP^2\right)\myDerA[\myP,\varyi]\myDerA[\myP,\varyj]
	+\frac{\myP}{2}\myDerA[\myP,\varyi]\myDerA[\myP,\varyj],
	\nonumber
\end{eqnarray}
and hence, from (\ref{eq::RealDXY}) and (\ref{eq::RealDXX}), we obtain
\begin{eqnarray}\label{eq::RelTwo}
	&&\myDerAA[\myP,\varzi,\varbzj]=
	\frac{1}{4}(1+\myP^2)\myDerAA[\myP,\varyi,\varyj]
	+\frac{\myP}{2}\myDerA[\myP,\varyi]\myDerA[\myP,\varyj].
\end{eqnarray}
Inserting (\ref{eq::RelOne}) and (\ref{eq::RelTwo})
in \ref{eq::DTwoArctan} we easily obtain
\begin{eqnarray}\nonumber
	&&\myDerAA[\arctan\myP,\varzi,\varbzj]=\frac{1}{4}\myDerAA[\myP,\varyi,\varyj],
	\nonumber
\end{eqnarray}
and this easily implies (\ref{eq::TubeLevi}).

It remains to prove that the function $\myP(z)$ is of class $C^2$ on
$\Cvx^{\rm ell}\setminus\Cvx$ and the equations
(\ref{eq::RealDX}), (\ref{eq::RealDXY}) and (\ref{eq::RealDXX}) hold.

Let $\auxF\in C^2(\R^n)$ be a global defining function for $\Cvx$,
that is $x\in\Cvx$ if, and only if, $\auxF(x)<0$ and
$(\myDerC[1]\auxF(x),\ldots,\myDerC[n]\auxF(x))\neq0$ if $x\in\partial\Cvx$.
Here $\myDerC[i]\auxF(x)$ is the derivative of $\auxF(x)$
with respect to the variable $\varxi$.

Let $\DiniF:\Cvx\times(\R^n\setminus\{0\})\times]0,+\infty[$ defined by
\begin{equation}\label{eq::BigFDef}
	\DiniF(x,y,\myP)=\auxF\left(x+p^{-1}y\right).
\end{equation}
We denote the derivatives $\myDerA[F,\varxi]$,
$\myDerA[F,\varyi]$ and $\myDerA[F,\myP]$ respectively as
$\myDerB[F,\varxi]$, $\myDerB[F,\varyi]$ and $\myDerB[F,\myP]$
(and similarilry for the higher order derivatives).

Since $\Cvx$ is convex then
$$
	\myDerB[\DiniF,\myP](x,y,\myP)=
	-\myP^{-2}\mySum y_\myIndexA\myDerC[\myIndexA]\auxF(x+p^{-1}y)
	\neq0
$$
when $x+p^{-1}y\in\partial\Cvx$ that is if $\myDerB[F,\myP](x,y,\myP)=0$.
Setting $z=x+iy$ the function $p(z)$ is characterized by the condition
\begin{equation}\label{eq::PDegZero}
	\DiniF(x,y,\myP(z))=0.
\end{equation}
By the Dini implicit function theorem, $p(z)$ is of class $C^2$
on $\DiniF:\Cvx\times(\R^n\setminus\{0\})$.

Taking the derivatives in (\ref{eq::PDegZero}) we obtain
\begin{eqnarray}
	\label{eq::PDegOneAA}
	&&\myDerB[\DiniF,\varxi]+\myDerB[\DiniF,\myP]\myDerB[\myP,\varxi]=0,
	\\
	\label{eq::PDegOneAB}
	&&\myDerB[\DiniF,\varyi]+\myDerB[\DiniF,\myP]\myDerB[\myP,\varyi]=0.
\end{eqnarray}
From (\ref{eq::BigFDef}) we obtain
\begin{eqnarray}
	\label{eq::PDegOneBA}
	&&\myDerB[\DiniF,\varxi]=\myDerC[i]\auxF,
	\\
	\label{eq::PDegOneBB}
	&&\myDerB[\DiniF,\varyi]=\myP^{-1}\myDerC[i]\auxF,
\end{eqnarray}
and hence
$$
	\myDerB[\myP,\varxi]=-\myDerB[\DiniF,\myP]^{-1}\myDerB[\DiniF,\varxi]
		=-\myP\myDerB[\DiniF,\myP]^{-1}\myDerB[\DiniF,\varyi]=\myP\myDerB[\myP,\varyi],
$$
which proves (\ref{eq::RealDX}).

Differentiating (\ref{eq::PDegOneAA}) with respect to $\varxj$ we compute
$$
	\myDerBB[\DiniF,\varxi,\varxj]
	   +\myDerBB[\DiniF,\varxi,\myP]\myDerB[\myP,\varxj]
	   +\myDerBB[\DiniF,\varxj,\myP]\myDerB[\myP,\varxi]
	   +\myDerBB[\DiniF,\myP,\myP]\myDerB[\myP,\varxi]\myDerB[\myP,\varxj]
	   +\myDerB[\DiniF,\myP]\myDerBB[\myP,\varxi,\varxj]=0,
$$
obtaining
\begin{equation}\label{eq::PDegwoAA}
	\myDerBB[\myP,\varxi,\varxj]=-\myDerB[\DiniF,\myP]^{-1}\left(
		\myDerBB[\DiniF,\varxi,\varxj]
	  	+\myDerBB[\DiniF,\varxi,\myP]\myDerB[\myP,\varxj]
	   	+\myDerBB[\DiniF,\varxj,\myP]\myDerB[\myP,\varxi]
	   	+\myDerBB[\DiniF,\myP,\myP]\myDerB[\myP,\varxi]\myDerB[\myP,\varxj]
	\right).
\end{equation}

Similarly we have
\begin{equation}\label{eq::PDegwoAB}
	\myDerBB[\myP,\varxi,\varyj]=-\myDerB[\DiniF,\myP]^{-1}\left(
		\myDerBB[\DiniF,\varxi,\varyj]
	  	+\myDerBB[\DiniF,\varxi,\myP]\myDerB[\myP,\varyj]
	   	+\myDerBB[\DiniF,\varyj,\myP]\myDerB[\myP,\varxi]
	   	+\myDerBB[\DiniF,\myP,\myP]\myDerB[\myP,\varxi]\myDerB[\myP,\varyj]
	\right)
\end{equation}
and
\begin{equation}\label{eq::PDegTwoAC}
	\myDerBB[\myP,\varyi,\varyj]=-\myDerB[\DiniF,\myP]^{-1}\left(
		\myDerBB[\DiniF,\varyi,\varyj]
	  	+\myDerBB[\DiniF,\varyi,\myP]\myDerB[\myP,\varyj]
	   	+\myDerBB[\DiniF,\varyj,\myP]\myDerB[\myP,\varyi]
	   	+\myDerBB[\DiniF,\myP,\myP]\myDerB[\myP,\varyi]\myDerB[\myP,\varyj]
	\right).
\end{equation}

But from (\ref{eq::PDegOneAA}) and (\ref{eq::PDegOneAB}) we have
\begin{eqnarray}
	&&\myDerBB[\DiniF,\varxi,\varxj]=\myDerC[i]\myDerC[j]\auxF,
	\nonumber\\
	&&\myDerBB[\DiniF,\varxi,\varyj]=\myP^{-1}\myDerC[i]\myDerC[j]\auxF,
	\nonumber\\
	&&\myDerBB[\DiniF,\varyi,\varyj]=\myP^{-2}\myDerC[i]\myDerC[j]\auxF,
	\nonumber\\
	&&\myDerBB[\DiniF,\varxi,\myP]=\myP^{-2}\mySum y_\myIndexA\myDerC[\myIndexA]\myDerC[i]\auxF,
	\nonumber\\
	&&\myDerBB[\DiniF,\varyi,\myP]=-\myP^{-2}\myDerC[i]\auxF
		-\myP^{-3}\mySum y_\myIndexA\myDerC[\myIndexA]\myDerC[i]\auxF,
	\nonumber
\end{eqnarray}
and hence
\begin{eqnarray}
	&&\myDerBB[\DiniF,\varxi,\varxj]=\myP^2\myDerBB[\DiniF,\varyi,\varyj],
	\nonumber\\
	&&\myDerBB[\DiniF,\varxi,\varyj]=\myP\myDerBB[\DiniF,\varyi,\varyj],
	\nonumber\\
	&&\myDerBB[\DiniF,\varxi,\myP]=\myP\myDerBB[\DiniF,\varyi,\varyj]+\myDerB[\DiniF,\varyi].
	\nonumber
\end{eqnarray}
Inserting such values in \ref{eq::PDegwoAB} we obtain
\begin{eqnarray}
	\myDerBB[\myP,\varxi,\varyj]
	&=&-\myDerB[\DiniF,\myP]^{-1}\left(
		\myDerBB[\DiniF,\varxi,\varyj]
	  	+\myDerBB[\DiniF,\varxi,\myP]\myDerB[\myP,\varyj]
	   	+\myDerBB[\DiniF,\varyj,\myP]\myDerB[\myP,\varxi]
	   	+\myDerBB[\DiniF,\myP,\myP]\myDerB[\myP,\varxi]\myDerB[\myP,\varyj]
	\right)
	\nonumber\\
	&=&-\myDerB[\DiniF,\myP]^{-1}\left(
		\myP\myDerBB[\DiniF,\varyi,\varyj]
	  	+(\myP\myDerBB[\DiniF,\varyi,\myP]+\myDerB[\DiniF,\varyi])\myDerB[\myP,\varyj]
	   	+\myP\myDerBB[\DiniF,\varyj,\myP]\myDerB[\myP,\varyi]
	   	+\myP\myDerBB[\DiniF,\myP,\myP]\myDerB[\myP,\varyi]\myDerB[\myP,\varyj]
	\right)
	\nonumber\\
	&=&-\myP\myDerB[\DiniF,\myP]^{-1}\left(
		\myDerBB[\DiniF,\varyi,\varyj]
	  	+\myDerBB[\DiniF,\varyi,\myP]\myDerB[\myP,\varyj]
	   	+\myDerBB[\DiniF,\varyj,\myP]\myDerB[\myP,\varyi]
	   	+\myDerBB[\DiniF,\myP,\myP]\myDerB[\myP,\varyi]\myDerB[\myP,\varyj]
	\right)-\myDerB[\DiniF,\myP]^{-1}\myDerB[\DiniF,\varyj]\myDerB[\myP,\varyj]
	\nonumber\\
	&=&\myP\myDerBB[\myP,\varyi,\varyi]+\myDerB[\myP,\varyj]\myDerB[\myP,\varyj],
	\nonumber
\end{eqnarray}
and this proves (\ref{eq::RealDXY}).

Finally, from (\ref{eq::PDegTwoAC})
\begin{eqnarray}
	\myDerBB[\myP,\varxi,\varxj]
	&=&-\myDerB[\DiniF,\myP]^{-1}\left(
		\myDerBB[\DiniF,\varxi,\varxj]
	  	+\myDerBB[\DiniF,\varxi,\myP]\myDerB[\myP,\varxj]
	   	+\myDerBB[\DiniF,\varxj,\myP]\myDerB[\myP,\varxi]
	   	+\myDerBB[\DiniF,\myP,\myP]\myDerB[\myP,\varxi]\myDerB[\myP,\varxj]
	\right)
	\nonumber\\
	&=&-\myDerB[\DiniF,\myP]^{-1}\left(
		\myP^2\myDerBB[\DiniF,\varyi,\varyj]
	  	+\myP(\myP\myDerBB[\DiniF,\varyi,\myP]+\myDerB[\DiniF,\varyi])\myDerB[\myP,\varyj]
	  \right.
	\nonumber\\
	&&\quad\quad\quad
		\left.
	  	+\myP(\myP\myDerBB[\DiniF,\varyj,\myP]+\myDerB[\DiniF,\varyj])\myDerB[\myP,\varyi]
	   	+\myP^2\myDerBB[\DiniF,\myP,\myP]\myDerB[\myP,\varyi]\myDerB[\myP,\varyj]
	\right)
	\nonumber\\
	&=&-\myP^2\myDerB[\DiniF,\myP]^{-1}\left(
		\myDerBB[\DiniF,\varyi,\varyj]
	  	+\myDerBB[\DiniF,\varyi,\myP]\myDerB[\myP,\varyj]
	   	+\myDerBB[\DiniF,\varyj,\myP]\myDerB[\myP,\varyi]
	   	+\myDerBB[\DiniF,\myP,\myP]\myDerB[\myP,\varyi]\myDerB[\myP,\varyj]
	\right)
	\nonumber\\
	&&\quad\quad\quad
	-\myP\myDerB[\DiniF,\myP]^{-1}\myDerB[\DiniF,\varyj]\myDerB[\myP,\varyj]
	-\myP\myDerB[\DiniF,\myP]^{-1}\myDerB[\DiniF,\varyi]\myDerB[\myP,\varyj]
	\nonumber\\
	&=&\myP^2\myDerBB[\myP,\varyi,\varyj]+2\myP\myDerB[\myP,\varyi]\myDerB[\myP,\varyj],
	\nonumber
\end{eqnarray}
obtaining hence (\ref{eq::RealDXX}).

The proof of the theorem is completed.
\qed



\end{document}